\documentclass[preprint,11pt]{elsarticle}
\usepackage{verbatim,a4margins,amsmath,bm,hyperref}
\journal{Journal of Computational and Applied Mathematics}

\newcommand{\dd}[2]{\frac{\mathrm d #1}{\mathrm d #2}}

\newcommand{\pdd}[2]{\frac{\partial #1}{\partial #2}}

\newcommand{\rem}[1]{}

\newcommand{\0}{\mathbf 0}

\newcommand{\av}{\mathbf a}
\newcommand{\f}{\mathbf f}
\newcommand{\uv}{\mathbf u}
\newcommand{\vv}{\mathbf v}
\newcommand{\w}{\mathbf w}
\newcommand{\x}{\mathbf x}

\newcommand{\cv}{\mathbf c}

\begin{document}

\begin{frontmatter}

\title{Higher-Order Corrections to\\ Optimisers based on Newton's Method}
\author[SJB]{Stephen Brooks}
\ead{sbrooks@bnl.gov}
\affiliation[SJB]{organization={Collider--Accelerator Department},
            addressline={Brookhaven National Laboratory}, 
            city={Upton},
            state={NY},
            postcode={11973}, 
            country={USA}}

\begin{abstract}
The Newton, Gauss--Newton and Levenberg--Marquardt methods all use the first derivative of a vector function (the Jacobian) to minimise its sum of squares.  When the Jacobian matrix is ill-conditioned, the function varies much faster in some directions than others and the space of possible improvement in sum of squares becomes a long narrow ellipsoid in the linear model.  This means that even a small amount of nonlinearity in the problem parameters can cause a proposed point far down the long axis of the ellipsoid to fall outside of the actual curved valley of improved values, even though it is quite nearby.  This paper presents a differential equation that `follows' these valleys, based on the technique of geodesic acceleration, which itself provides a 2$^\mathrm{nd}$ order improvement to the Levenberg--Marquardt iteration step.  Higher derivatives of this equation are computed that allow $n^\mathrm{th}$ order improvements to the optimisation methods to be derived.  These higher-order accelerated methods up to 4$^\mathrm{th}$ order are tested numerically and shown to provide substantial reduction of both number of steps and computation time.
\end{abstract}

\end{frontmatter}

\section{Definitions and Introduction}
Consider finding the value of a vector $\x$ such that the vector-valued function $\f(\x)=\0$, noting the input and output of $\f$ might have different dimensions.

Newton's method solves $J(\x)(\x_\mathrm{new}-\x)=-\f(\x)$ where $J$ is the Jacobian matrix.  The Gauss--Newton algorithm generalises this to rectangular $J$ using pseudo-inverses that may be calculated using Singular Value Decomposition (SVD).  The Levenberg--Marquardt algorithm \cite{Levenberg,Marquardt} introduces a damping factor into this pseudo-inverse, which allows progress along `easier' directions without having to go far in `difficult' directions that may exhibit nonlinearity.

The remainder of this paper will be written for the simpler Newton method, where $J^{-1}$ is the inverse of the square Jacobian matrix.  However, the algorithms derived also work for the Gauss-Newton pseudo-inverse $[J^{-1}]_{GN}=(J^TJ)^{-1}J^T$ and the damped Levenberg--Marquardt version $[J^{-1}]_{LM(\lambda)}=(J^TJ+\lambda I)^{-1}J^T$.  The latter is used in numerical tests.

One common source of slow convergence is that $J$ is ill-conditioned, so the optimisation valley is much narrower in some directions than others, while the problem contains some nonlinearity, which may seem small but is amplified once you change into coordinates where $J$ is well-conditioned.  This is because even a small amount of nonlinearity can make the long narrow valley stop overlapping with its approximation in the linear model.  This reduced range of validity of the linear model means many small steps have to be taken.

\section{Natural Optimisation Pathway}
The goal of the Newton step is to reduce the error vector $\f$, ideally to zero.  For a nonlinear function, the optimisation follows a curved pathway \cite{GeodesicAccelLett,GeodesicAccelE} and one natural such pathway is $\x(t)$ defined implicitly by
\[ \f(\x(t)) = (1-t)\f(\x(0)) \]
for $t\in[0,1]$.  This scales down all components of the error equally and at $t=1$ it reaches the true solution.

Taking the first derivative of this equation gives
\[ \sum_i\partial_i\f(\x(t))\dot x_i(t) = J(\x(t))\dot\x(t) = -\f(\x(0)) = -\frac{\f(\x(t))}{1-t}, \]
which at $t=0$ makes $\dot\x$ equal to the Newton step and to a scaling of it for all $0<t<1$.  So this pathway is always tangent to the Newton step direction and corresponds to the limit of a Newton algorithm run with steps scaled down to be infinitesimally small.

\section{Higher-Order Derivatives}
If the pathway curves, one may wonder if longer steps can be taken if the curvature is taken into account.  The second and higher derivatives of the equation defining the natural pathway have the form
\[ \dd{^n}{t^n}\f(\x(t)) = \0 \]
for $n\ge 2$.  Multiple derivatives of a function composition ($\f\circ\x$ here) are given by Faà di Bruno's formula \cite{Abrogast,FaadiBruno1,FaadiBruno2}
\[ \dd{^n}{t^n}\f(\x(t)) = \sum_{\pi\in\Pi_n}\f^{(|\pi|)}(\x(t))\bigotimes_{p\in\pi}\x^{(|p|)}(t), \]
where $\Pi_n$ is the set of all partitions of $\{1,2,...,n\}$.  The $d^\mathrm{th}$ derivative of the vector function $\f$ is a tensor that takes $d$ vectors as input and outputs a vector, with elements defined by
\[ f^{(d)}(\x)^i_{j_1j_2...j_d} = \pdd{^df_i(\x)}{x_{j_1} \partial x_{j_2} ... \partial x_{j_d}}. \]
Note that $\f^{(1)}=J$.  This paper will adopt compact notation where tensor products of vectors $\uv\otimes\vv\otimes\w$ will be written $\uv\vv\w$ so that $(\uv\vv\w)_{ijk}=u_iv_jw_k$.  These may be contracted with the derivative tensor to give a vector written in the form $\f^{(3)}\uv\vv\w$, where $(\f^{(3)}\uv\vv\w)_n=\sum_{i,j,k}f^{(3)n}_{ijk}u_iv_jw_k$.
 
\subsection{Second Order}
For $n=2$, $\Pi_2=\{\{\{1\},\{2\}\},\{\{1,2\}\}\}$ and 
\[ \dd{^2}{t^2}\f(\x(t)) = \f^{(2)}(\x(t))\dot\x(t)\dot\x(t) + \f^{(1)}(\x(t))\ddot\x(t) = \0 \]
\[ \Rightarrow \qquad \ddot\x(t) = -J^{-1}(\x(t))\f^{(2)}(\x(t))\dot\x(t)\dot\x(t). \]
This agrees with the well-known \cite{GeodesicAccelLett,GeodesicAccelE,GeodesicArxiv1,GeodesicArxiv2} quadratic acceleration term for Levenberg--Marquardt if the $J^{-1}$ is replaced by a damped pseudo-inverse $[J^{-1}]_{LM(\lambda)}$.

 \subsection{Third Order}
 For conciseness, $\x$ and its derivatives will be evaluated at $t=0$ unless otherwise stated and $\f$ and its derivatives at $\x$.  For $n=3$, 
 \[ \Pi_3=\{\{\{1\},\{2\},\{3\}\},\{\{1\},\{2,3\}\},\{\{2\},\{1,3\}\},\{\{3\},\{1,2\}\},\{\{1,2,3\}\}\} \] and 
\[ \dd{^3}{t^3}\f = \f^{(3)}\dot\x\dot\x\dot\x + 3\f^{(2)}\dot\x\ddot\x + \f^{(1)}\x^{(3)} = \0. \]
This gives the third derivative of $\x$ as 
\[ \x^{(3)} = -J^{-1}( \f^{(3)}\dot\x\dot\x\dot\x + 3\f^{(2)}\dot\x\ddot\x). \]
 
\subsection{Fourth Order, Recurrence and General Case}
Higher-order expressions can be obtained either from the set partitions $\Pi_n$ or the equivalent differentiation chain and product rules that obtain $\dd{^{n+1}}{t^{n+1}}\f$ from $\dd{^n}{t^n}\f$.  The formulae
\[ \dd{}{t}\f^{(n)} = \f^{(n+1)}\dot\x \qquad \mathrm{and} \qquad \dd{}{t}\x^{(n)} = \x^{(n+1)} \]
together with the product rule are enough to generate the full sequence.  Starting from $n=2$,
\begin{eqnarray*}
 \f^{(2)}\dot\x\dot\x + \f^{(1)}\ddot\x &=& \0 \\
 \f^{(3)}\dot\x\dot\x\dot\x + 3\f^{(2)}\dot\x\ddot\x + \f^{(1)}\x^{(3)} &=& \0 \\
 \f^{(4)}\dot\x\dot\x\dot\x\dot\x + 6\f^{(3)}\dot\x\dot\x\ddot\x + 4\f^{(2)}\dot\x\x^{(3)} + 3\f^{(2)}\ddot\x\ddot\x + \f^{(1)}\x^{(4)} &=& \0
\end{eqnarray*}
and so on.  A computer algebra system can generate these terms based on a rule like
\[ \dd{}{t}\f^{(n)}\x^{(a)}\x^{(b)}\x^{(c)} =  \]
\[ \f^{(n+1)}\x^{(1)}\x^{(a)}\x^{(b)}\x^{(c)} +
\f^{(n)}\x^{(a+1)}\x^{(b)}\x^{(c)} +
\f^{(n)}\x^{(a)}\x^{(b+1)}\x^{(c)} +
\f^{(n)}\x^{(a)}\x^{(b)}\x^{(c+1)} \]
and collecting like terms, for example by sorting the $\x$ derivatives in increasing order.
 
The highest derivative $\x^{(n)}$ may be moved to the other side to get a formula like
\[ \x^{(4)} = -J^{-1}( \f^{(4)}\dot\x\dot\x\dot\x\dot\x + 6\f^{(3)}\dot\x\dot\x\ddot\x + 4\f^{(2)}\dot\x\x^{(3)} + 3\f^{(2)}\ddot\x\ddot\x ), \]
shown for the $n=4$ case, which expresses it in terms of lower derivatives of $\x$.

\section{Taking Finite Steps}
The derivatives $\x^{(n)}$ calculated above can produce a corrected higher-order step using the Taylor series of $\x$ around $t=0$
\[ \x(\epsilon) = \sum_{n=0}^\infty \frac{1}{n!}\epsilon^n\x^{(n)}, \]
where the step is thought of as stopping at a time $t=\epsilon$ in the parameterisation of the natural pathway.  This unknown $\epsilon$ may seem like a problem but it can be made to cancel.  Define the correction at order $n$ to be the $n^\mathrm{th}$ term of the Taylor series:
\[ \cv_n=\frac1{n!}\epsilon^n\x^{(n)}. \]
The step begins at $\cv_0=\x$ and the first order uncorrected step ends at $\cv_0+\cv_1$, so has length $\cv_1$.
Now recall that for $n\ge 2$, the derivatives of $\f\circ\x$ are zero and use Faà di Bruno's formula as before:
\[ \dd{^n}{t^n}\f(\x(t)) = \sum_{\pi\in\Pi_n}\f^{(|\pi|)}(\x(t))\bigotimes_{p\in\pi}\x^{(|p|)}(t) = \0. \]
Multiplying both sides by $\epsilon^n$ gives
\[ \sum_{\pi\in\Pi_n}\f^{(|\pi|)}(\x(t))\bigotimes_{p\in\pi}\epsilon^{|p|}\x^{(|p|)}(t) = \0, \]
using the fact that $\pi$ is a partition of $\{1,2,...,n\}$, so the sum of sizes $|p|$ of all its elements is $n$.  Noting that $\epsilon^n\x^{(n)}=n!\cv_n$ and evaluating at $t=0$ gives
\[ \sum_{\pi\in\Pi_n}\f^{(|\pi|)}\bigotimes_{p\in\pi}|p|!\cv_{|p|} = \0. \]
This formula is the basis for calculating corrections $\cv_n$ for finite steps in the following sections.

\subsection{The Meaning of $\epsilon$}
Observant readers might have noticed that $\cv_1=\epsilon\dot\x$ and in an earlier section, $\dot\x=-J^{-1}\f$, so taking a full Newton step would imply $\epsilon=1$.  This paper treats $\epsilon$ as a small value because when experiencing slow convergence from the `narrow curving valleys' problem, the area of validity for the local linear model (the trust region) is much smaller than what is required to go all the way to the model minimum.  This means the steps taken that succeed in reducing the function sum of squares would only be a fraction of the Newton step, for example a Levenberg--Marquardt step with $\lambda$ chosen large enough to damp away the longest-range movement axes of the exact Newton scheme.

\section{Finite Difference Schemes}
The higher-order corrections $\cv_n$ are expressible in terms of multiple directional derivatives of $\f$.  For a numerical method, these derivatives must be calculated from function values, or at most, the Jacobian used by the algorithm.  In this paper finite difference schemes are used, some of which have their `stencils' of sampled points spread in multiple axes to give mixed derivatives.

It would also be possible to use an automatic differentiation scheme here, provided it supports higher order and mixed derivatives.  However, because each $\cv_i$ is $O(\epsilon^i)$ by definition and is never divided by $\epsilon$, finite difference schemes are not expected to result in a dramatic loss of accuracy in the following algorithms.

\subsection{Second Order}
For $n=2$, the general formula gives
\[ \f^{(2)}\cv_1\cv_1 + \f^{(1)}2\cv_2 = \0 \]
\[ \Rightarrow \qquad \cv_2 = -\tfrac12J^{-1}\f^{(2)}\cv_1\cv_1. \]
Taylor expansion of $\f$ in the direction $\cv_1$ of the original uncorrected step gives
\[ \f(\x+\cv_1) = \f + J\cv_1 + \tfrac12\f^{(2)}\cv_1\cv_1 + O(\epsilon^3) \]
\[ \Rightarrow \qquad \tfrac12\f^{(2)}\cv_1\cv_1 = \f(\x+\cv_1) - (\f + J\cv_1) + O(\epsilon^3). \]
In other words, the difference between $\f(\x+\cv_1)$ and a linear estimate using the $\f$ and $J$ already calculated at $\x(0)$, is to leading order a second derivative term similar to the one required for calculating $\cv_2$.  Thus,
\[ \cv_2 = -J^{-1}(\f(\x+\cv_1) - (\f + J\cv_1)) + O(\epsilon^3). \]

\begin{figure}[!htb]
   \centering
   \includegraphics*[width=0.5\columnwidth]{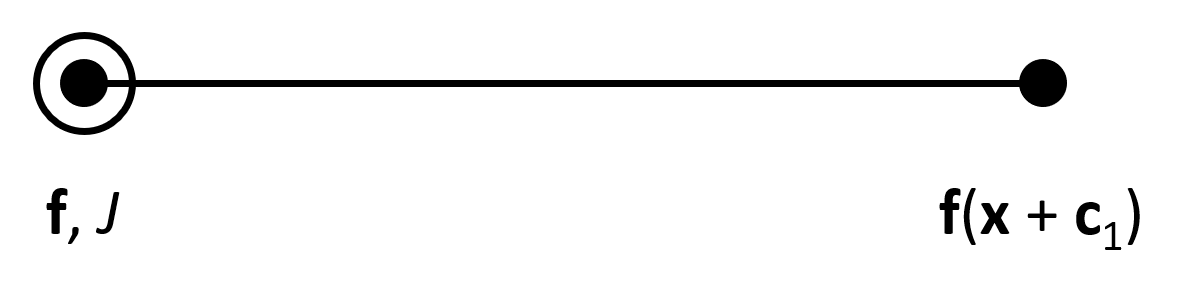}
   \caption{Finite difference stencil for calculating the second order correction $\cv_2$.  Points represent evaluations of the function $\f$ and rings represent evaluations of its Jacobian.}
   \label{fig:stencil2}
\end{figure}

The evaluations required for this calculation are shown in Figure \ref{fig:stencil2}.  In this case, only one other point besides the evaluations of $\f$ and $J$ at $\x$ is needed.

\subsection{Third Order}
For $n=3$, the general formula gives
\[ \f^{(3)}\cv_1\cv_1\cv_1 + 3\f^{(2)}\cv_12\cv_2 + \f^{(1)}6\cv_3 = \0. \]
\[ \Rightarrow \qquad \cv_3 = -\tfrac16J^{-1}(\f^{(3)}\cv_1\cv_1\cv_1 + 6\f^{(2)}\cv_1\cv_2). \]
There are a few differences from the second order case:
\begin{itemize}
\setlength{\itemsep}{0pt}
\item There is a third order derivative $\f^{(3)}\cv_1\cv_1\cv_1$, which will require an additional stencil point in the direction of $\cv_1$.
\item Errors will now have to be $O(\epsilon^4)$ as the main terms have size $O(\epsilon^3)$.
\item There is a mixed derivative $\f^{(2)}\cv_1\cv_2$, requiring a two dimensional stencil pattern.
\end{itemize}
The mixed derivative requires knowledge of the direction $\cv_2$, which must be evaluated first.  The second order stencil for $\cv_2$ had error $O(\epsilon^3)$ and now $O(\epsilon^4)$ is needed, so even the lower-order derivative $\f^{(2)}\cv_1\cv_1$ will have to be evaluated using a third order stencil.  Fortunately, this stencil is also needed for evaluating $\f^{(3)}\cv_1\cv_1\cv_1$, so all coefficients of a cubic approximation to $\f$ in this direction can be known.

\subsubsection{Phase One: Calculating $\cv_2$}
The additional stencil point in the $\cv_1$ direction is chosen to be $\x+\frac12\cv_1$ here, although other choices are possible.  To third order,
\[ \f(\x+\tfrac12\cv_1) = \f + \tfrac12J\cv_1 + \tfrac18\f^{(2)}\cv_1\cv_1 + \tfrac1{48}\f^{(3)}\cv_1\cv_1\cv_1 + O(\epsilon^4) \]
\[ \f(\x+\cv_1) = \f + J\cv_1 + \tfrac12\f^{(2)}\cv_1\cv_1 + \tfrac16\f^{(3)}\cv_1\cv_1\cv_1 + O(\epsilon^4). \]
Writing the nonlinear part of $\f$ as $\f_{nl}(\x+\av)=\f(\x+\av) - (\f + J\av)$, the derivatives in the $\cv_1$ direction can be expressed
\[ \f^{(2)}\cv_1\cv_1 = 16\f_{nl}(\x+\tfrac12\cv_1) - 2\f_{nl}(\x+\cv_1) + O(\epsilon^4) \]
\[ \f^{(3)}\cv_1\cv_1\cv_1 = 12\f_{nl}(\x+\cv_1) - 48\f_{nl}(\x+\tfrac12\cv_1) + O(\epsilon^4) \]
and $\cv_2$ calculated from the formula
$\cv_2 = -\frac12J^{-1}\f^{(2)}\cv_1\cv_1$
with $O(\epsilon^4)$ error.

\subsubsection{Phase Two: Calculating $\cv_3$}
This step requires the mixed derivative $\f^{(2)}\cv_1\cv_2$.  Expressions of the form $\f^{(3)}\uv\vv\w=(\uv\cdot\nabla)(\vv\cdot\nabla)(\w\cdot\nabla)\f$ are iterated directional derivatives.  Each directional derivative can be approximated to leading order as
\begin{eqnarray*}
(\uv\cdot\nabla)\f(\x) &=& \frac{\f(\x+\epsilon\uv)-\f(\x)}{\epsilon}+O(\epsilon) \\
\Rightarrow \qquad (\epsilon\uv\cdot\nabla)\f(\x) &=& \f(\x+\epsilon\uv)-\f(\x)+O(\epsilon^2) \\
\Rightarrow \qquad (\epsilon^n\uv\cdot\nabla)\f(\x) &=& \f(\x+\epsilon^n\uv)-\f(\x)+O(\epsilon^{2n}).
\end{eqnarray*}
In the last formula above, $\epsilon^n\uv$ represents an $O(\epsilon^n)$ sized term such as $\cv_n$.  Using this multiple times allows mixed derivatives to be expressed to leading order as combinations of function evaluations at different points (i.e. finite difference stencils).  For example,
\begin{eqnarray*}
\f^{(2)}\cv_1\cv_2 &=& (\cv_1\cdot\nabla)(\cv_2\cdot\nabla)\f \\
&\simeq&  (\cv_1\cdot\nabla)(\f(\x+\cv_2) - \f(\x)) \\
&\simeq& \f(\x+\cv_2+\cv_1) - \f(\x+\cv_1) - (\f(\x+\cv_2) - \f(\x)).
\end{eqnarray*}
Here, all terms have size $O(\epsilon^3)$ and all approximations are leading order accurate meaning the error is no worse than $O(\epsilon^4)$, as required.

In general, a $d^\mathrm{th}$ derivative of different directions would require $2^d$ evaluations.  An $n$ times repeated derivative in the same direction only requires $n+1$ as some of the evaluation points are coincident.  A mixture like $\f^{(a+b+c)}\uv^{\otimes a}\vv^{\otimes b}\w^{\otimes c}$ would require evaluation at $(a+1)(b+1)(c+1)$ points.

Some efficiencies may be gained from coincident evaluation points and the fact that the full first derivative $J$ is usually evaluated at $\x$ already.  This was used in the previous `second order' section, which only required one additional evaluation point for $\f^{(2)}\cv_1\cv_1$ rather than three.

\begin{figure}[!htb]
   \centering
   \includegraphics*[width=0.75\columnwidth]{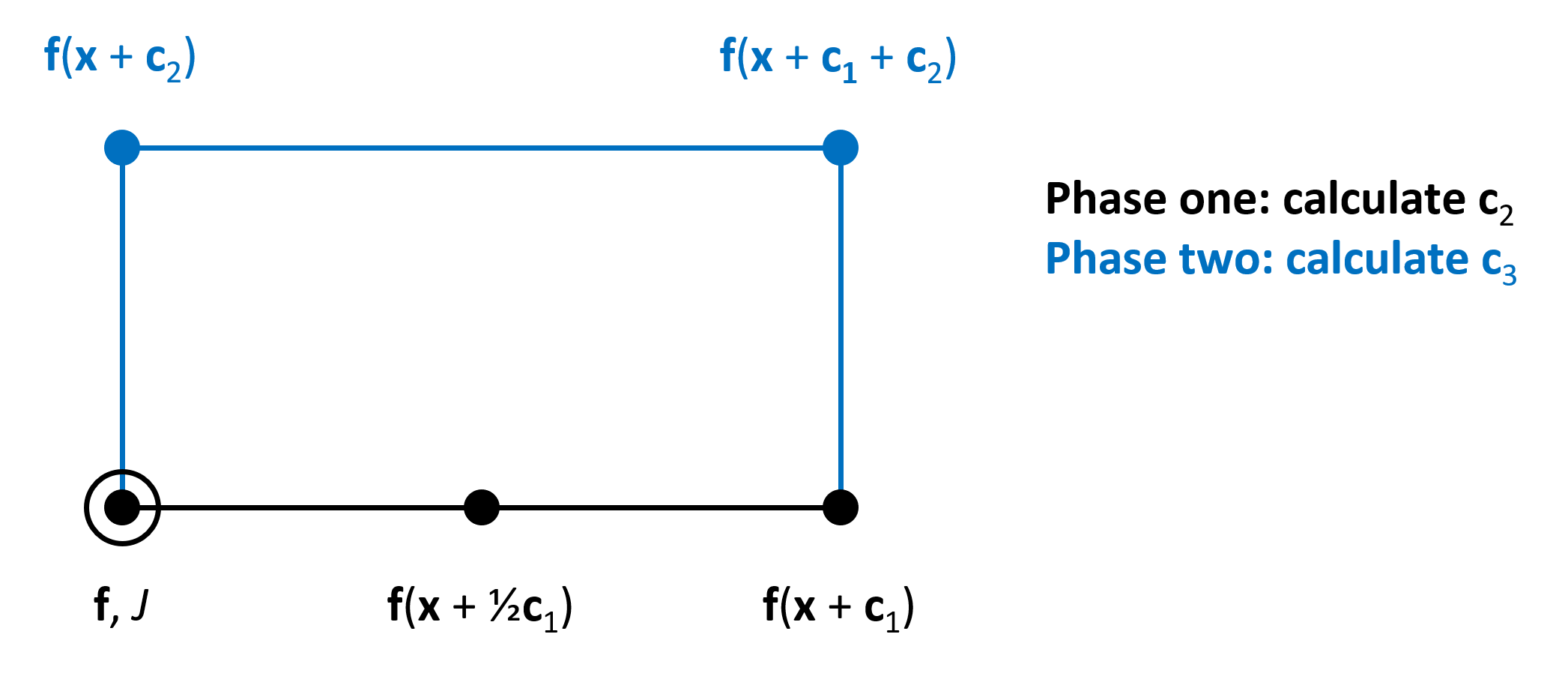}
   \caption{Finite difference stencil for calculating the third order correction $\cv_3$ along with the second order correction $\cv_2$ that is also required.}
   \label{fig:stencil3}
\end{figure}

Now everything is in place to evaluate $\cv_3 = -\frac16J^{-1}(\f^{(3)}\cv_1\cv_1\cv_1 + 6\f^{(2)}\cv_1\cv_2)$.  The full step will be corrected from $\cv_1$ to $\cv_1+\cv_2+\cv_3$.  The evaluations required are shown in Figure \ref{fig:stencil3}.

\subsection{Fourth Order}
For $n=4$, the general formula gives
\[ \f^{(4)}\cv_1\cv_1\cv_1\cv_1 + 6\f^{(3)}\cv_1\cv_12\cv_2 + 4\f^{(2)}\cv_16\cv_3 + 3\f^{(2)}2\cv_22\cv_2 + \f^{(1)}24\cv_4 = \0 \]
\[ \Rightarrow \qquad \cv_4 = -\tfrac1{24}J^{-1}(\f^{(4)}\cv_1\cv_1\cv_1\cv_1 + 12\f^{(3)}\cv_1\cv_1\cv_2 + 24\f^{(2)}\cv_1\cv_3 + 12\f^{(2)}\cv_2\cv_2). \]
As expected, there are more higher-order and mixed derivatives.  The double derivative $\f^{(2)}\cv_2\cv_2$ can take advantage of the Jacobian to eliminate a point from the stencil, just as previous unidirectional derivatives did.  The direction $\cv_3$ is now involved in the derivatives, so three evaluation phases are required.  All errors have to be $O(\epsilon^5)$ including those of $\cv_2$ and $\cv_3$.

\subsubsection{Phase One: Calculating $\cv_2$}
An additional stencil point $\x+\frac32\cv_1$ will be added to increase the order of accuracy in the $\cv_1$ direction.  Defining $\f_{nl}(\x+\av)=\f(\x+\av) - (\f + J\av)$ as before, 
\begin{eqnarray*}
\f^{(2)}\cv_1\cv_1 &\simeq& 24\f_{nl}(\x+\tfrac12\cv_1) - 6\f_{nl}(\x+\cv_1) + \tfrac89\f(\x+\tfrac32\cv_1) \\
\f^{(3)}\cv_1\cv_1\cv_1 &\simeq& -120\f_{nl}(\x+\tfrac12\cv_1) + 48\f_{nl}(\x+\cv_1) -8\f(\x+\tfrac32\cv_1) \\
\f^{(4)}\cv_1\cv_1\cv_1\cv_1 &\simeq& 192\f_{nl}(\x+\tfrac12\cv_1) -96\f_{nl}(\x+\cv_1) +\tfrac{64}3\f(\x+\tfrac32\cv_1),
\end{eqnarray*}
all with errors $O(\epsilon^5)$.  These formulae came from writing out the Taylor expansions and inverting the system of equations, which can also be done by inverting a matrix as the equations are linear in $\f$ and its derivatives.  Calculate $\cv_2 = -\frac12J^{-1}\f^{(2)}\cv_1\cv_1$ using the first formula above.

\subsubsection{Phase Two: Calculating $\cv_3$}
The mixed derivative $\f^{(3)}\cv_1\cv_1\cv_2$ needs a grid with three points in the $\cv_1$ direction and two in the $\cv_2$ direction for a total of six.  Many previous points can be re-used, with the only new point for fourth order in this phase being $\x+\frac12\cv_1+\cv_2$.
\begin{eqnarray*}
\f^{(3)}\cv_1\cv_1\cv_2 &\simeq& (4\f(\x+\cv_2)-8\f(\x+\tfrac12\cv_1+\cv_2)+4\f(\x+\cv_1+\cv_2)) - \\
	&& (4\f-8\f(\x+\tfrac12\cv_1)+4\f(\x+\cv_1)) \\
\f^{(2)}\cv_1\cv_2 &\simeq& (-3\f(\x+\cv_2)+4\f(\x+\tfrac12\cv_1+\cv_2)-\f(\x+\cv_1+\cv_2)) - \\
	&& (-3\f+4\f(\x+\tfrac12\cv_1)-\f(\x+\cv_1)) \\
\f^{(2)}\cv_2\cv_2 &\simeq& 2\f_{nl}(\x+\cv_2).
\end{eqnarray*}
Again, all errors are $O(\epsilon^5)$ and $\cv_3 = -\frac16J^{-1}(\f^{(3)}\cv_1\cv_1\cv_1 + 6\f^{(2)}\cv_1\cv_2)$ can now be calculated to the required accuracy.

\subsubsection{Phase Three: Calculating $\cv_4$}
\begin{figure}[!htb]
   \centering
   \includegraphics*[width=0.85\columnwidth]{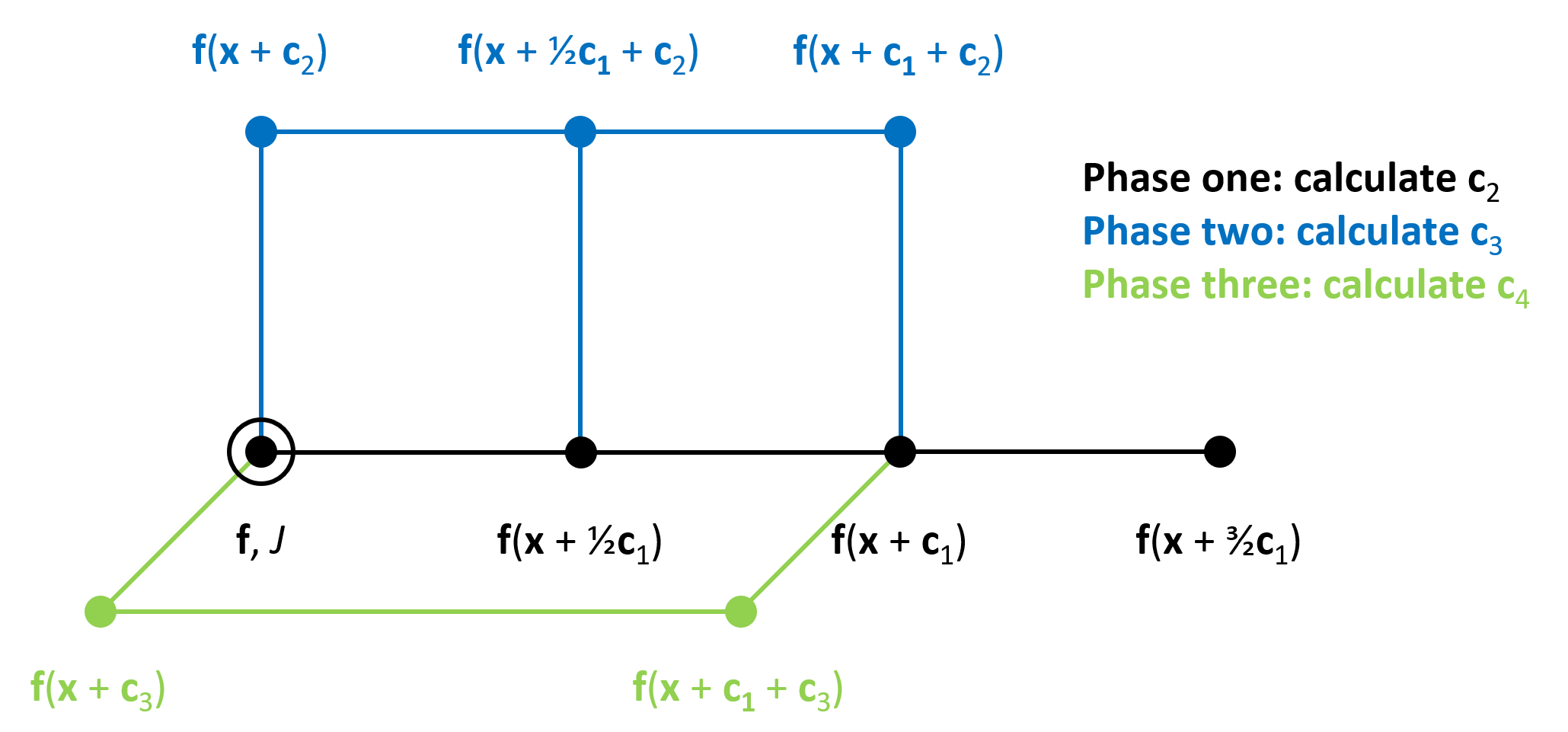}
   \caption{Finite difference stencil for calculating the second, third and fourth order corrections $\cv_2+\cv_3+\cv_4$ at fourth order accuracy.}
   \label{fig:stencil4}
\end{figure}

This phase requires the single mixed derivative $\f^{(2)}\cv_1\cv_3$, which can be handled analogously to when $\f^{(2)}\cv_1\cv_2$ first appeared, by extending the stencil in the $\cv_3$ direction with the two points $\x+\cv_3$ and $\x+\cv_1+\cv_3$.
\begin{eqnarray*}
\f^{(2)}\cv_1\cv_3 &\simeq& \f(\x+\cv_1+\cv_3) - \f(\x+\cv_3) - (\f(\x+\cv_1) - \f(\x)).
\end{eqnarray*}
Now all values are available to evaluate the fourth order correction $\cv_4 = -\tfrac1{24}J^{-1}(\f^{(4)}\cv_1\cv_1\cv_1\cv_1 + 12\f^{(3)}\cv_1\cv_1\cv_2 + 24\f^{(2)}\cv_1\cv_3 + 12\f^{(2)}\cv_2\cv_2)$.  The full step will be corrected from $\cv_1$ to $\cv_1+\cv_2+\cv_3+\cv_4$ and the evaluations required are shown in Figure \ref{fig:stencil4}.

It is clear that this process could be continued to even higher orders, although the stencils would require more and more points.  Practically, automated selection of points and calculations of the stencil coefficients would also be required.

\section{Numerical Test Problem}
The higher-order algorithms were tested on a simple function to verify their performance.  This function had to exhibit the `narrow curving valleys' problem in its sum of squares, so a very anisotropic function $(x,y)\mapsto(x,Ky)$ for $K\gg1$ was chosen and then some nonlinearity added in parameter space.  The resulting function is
\[ \f(x,y)=(x+y^2,K(y-x^2)). \]
Typical values of $K=10^6$ were used and the iteration was started from the arbitrary point $(x,y)=(\pi,e)$, moving towards the minimum sum of squares at $\f(0,0)=(0,0)$.  

\begin{figure}[!htb]
   \centering
   \includegraphics*[width=\columnwidth]{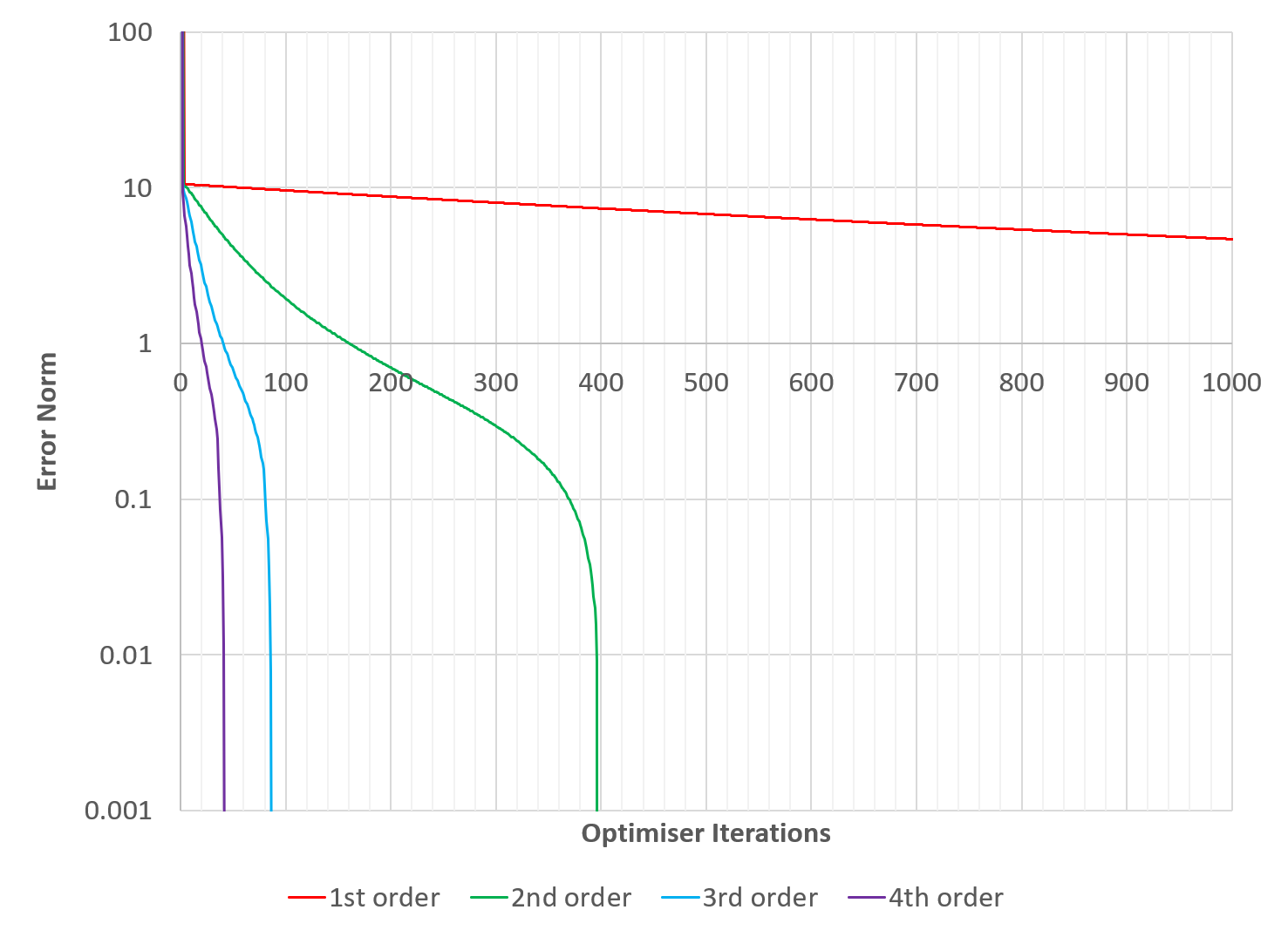}
   \caption{Performance of higher-order corrected Levenburg--Marquardt methods on a test problem.}
   \label{fig:testproblem}
\end{figure}

Figure \ref{fig:testproblem} shows the improved performance of the higher order corrected methods on the test problem with $K=10^6$.  The error norm in the plot is the value of $|\f(x,y)|$ after each iteration.  There is a slow convergence region for $0.1\le|\f|\le10$ after which the algorithm converges rapidly.  When the full Newton step using the linear model becomes a valid, convergence should be quadratic with the error norm roughly squaring on each iteration.  This rapid convergence appears as the near-vertical descending lines on the graph for $|\f|<0.1$.

Source code in C for this example is available at \cite{src}.

\subsection{Varying Valley Anisotropy}
Varying $K$ should show the relative performance of the different order algorithms as the valley gets narrower, while the curvature is kept constant.  Table \ref{tab:anisotropy} shows the number of iterations required to converge as $K$ is varied between $1$ and $10^{12}$.

\begin{table}[!hbt]
   \centering
   \caption{Test problem convergence times for different order methods as the anisotropy factor $K$ is varied.}
   \begin{tabular}{lcccc}
       \hline
\textbf{Anisotropy $K$} & \textbf{1st order} & \textbf{2nd order} & \textbf{3rd order} & \textbf{4th order} \\
       \hline
1 & 8 & 6 & 5 & 5 \\
10 & 15 & 8 & 6 & 5 \\
100 & 47 & 16 & 9 & 8 \\
1000 & 196 & 30 & 18 & 11 \\
10000 & 880 & 68 & 24 & 18 \\
100000 & 4041 & 162 & 50 & 27 \\
$10^{6}$ & 18733 & 397 & 88 & 43 \\
$10^{7}$ & $>$20000 & 971 & 166 & 70 \\
$10^{8}$ & $>$20000 & 2432 & 312 & 110 \\
$10^{9}$ & $>$20000 & 5828 & 631 & 243 \\
$10^{10}$ & $>$20000 & $>$20000 & 2876 & 968 \\
$10^{11}$ & $>$20000 & $>$20000 & 10886 & 2706 \\
$10^{12}$ & $>$20000 & $>$20000 & $>$20000 & 9159 \\
       \hline
   \end{tabular}
   \label{tab:anisotropy}
\end{table}

A simplified model is that the valley has width $1/K$ while the $n^\mathrm{th}$ order method has error term $O(\epsilon^{n+1})$, so this error would push the proposed step out of the valley when $O(\epsilon^{n+1})=1/K$ or $\epsilon=O(K^{-\frac{1}{n+1}})$.  This would give $O(1/\epsilon)=O(K^{\frac{1}{n+1}})$ steps to convergence.

\begin{figure}[!htb]
   \centering
   \includegraphics*[width=\columnwidth]{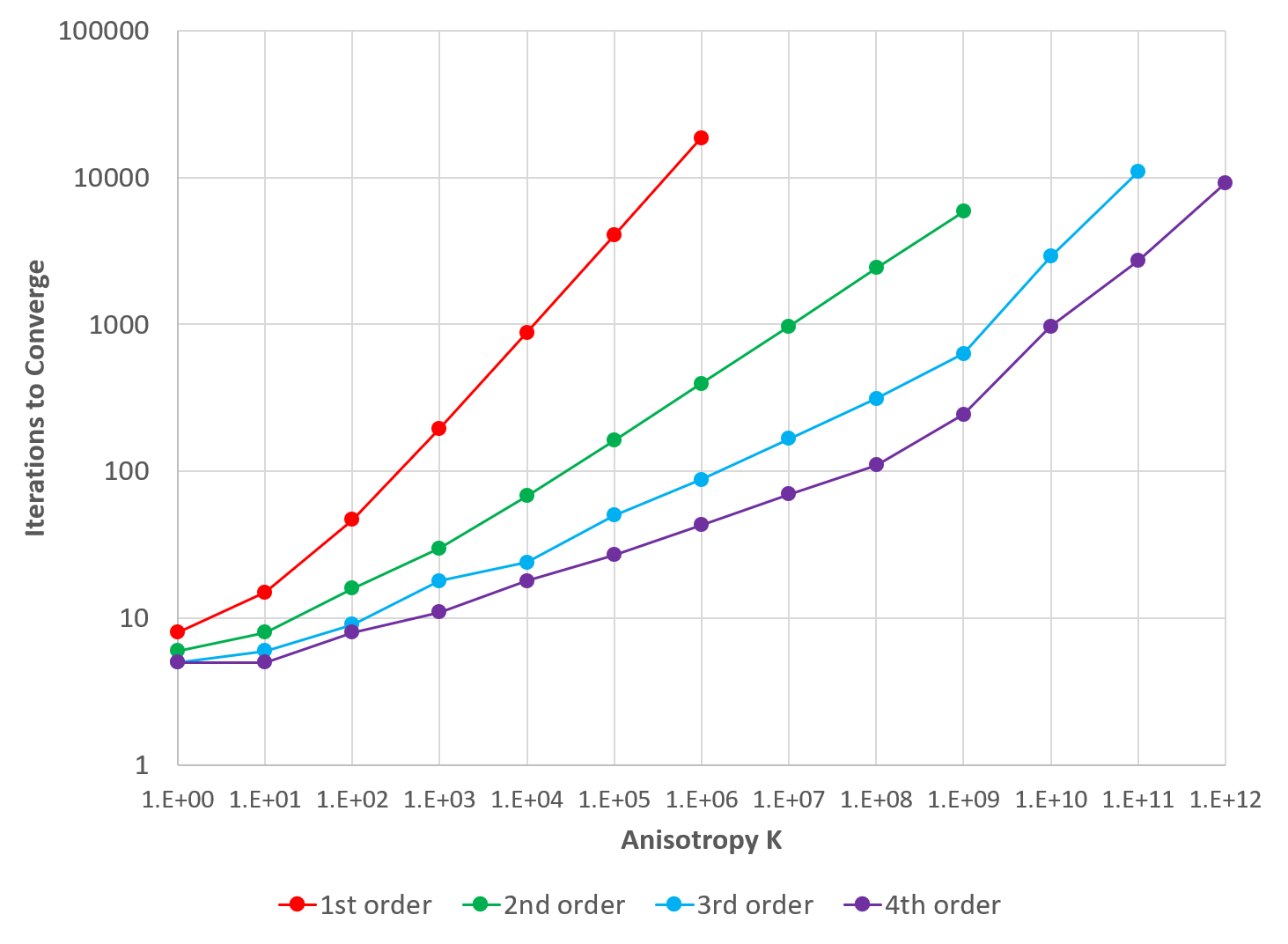}
   \caption{Test problem convergence times for different order methods as the anisotropy factor $K$ is varied.}
   \label{fig:anisotropy}
\end{figure}

Plotting the data on a log-log plot in Figure \ref{fig:anisotropy} reveals straight lines in parts of the data that suggest a power law relationship.  For $K\ge10^9$ there is an additional increase in convergence time, which may be from the limits of double precision used for the calculation.  Taking the gradient through the last three available points with $K\le10^8$ gives power law exponents of $0.660, 0.392, 0.265, 0.203$ for the first through fourth order methods.  This is somewhat similar to the simplified model's $\frac12, \frac13, \frac14, \frac15$ although lower-order methods seem slower, particularly the first order.

\clearpage
\subsection{Combination with Levenburg--Marquardt Method}
These numerical experiments were performed with a Levenburg--Marquardt method enhanced with the higher-order corrections.  Step length $\epsilon$ is controlled by choosing the damping parameter $\lambda\ge 0$ in the pseudo-inverse $[J^{-1}]_{LM(\lambda)}=(J^TJ+\lambda I)^{-1}J^T$.  The note \cite{SVDnote} shows that $\lambda=0$ gives the full Gauss--Newton step, while $\lambda\rightarrow\infty$ produces infinitesimal steepest gradient steps, with the values of $\lambda$ in between producing optimal reductions in the linear model for a given step size.

A good choice of $\lambda$ should be chosen at each step.  In this study the values 
\[ \lambda_n = \lambda_\mathrm{old} 10000^{(n/10)^3} \qquad \mbox{for} \quad -10\le n \le 10, \]
where $\lambda_\mathrm{old}$ is the value from the previous step, are run in parallel and the one that produces the best reduction in $|\f|$ chosen.  The initial step uses $\lambda_\mathrm{old}=1$.

The fact these steps are run in parallel means that for the higher order methods, many initial Levenburg--Marquardt steps $\cv_1(\lambda_n)$ are calculated, each of which has higher order corrections $\cv_i(\lambda_n)$.  The function $\f$ is evaluated at all the corrected step points $\x_\mathrm{out}(\lambda_n)=\sum_{i=0}^{order}\cv_i(\lambda_n)$ and the one with lowest $|\f|$ and its corresponding $\lambda_n$ is chosen.

The higher order methods enable longer steps and thus smaller values of $\lambda$ to be used.  The scheme above is somewhat wasteful by trying 21 values of $\lambda$ each step, but on modern computers these can be parallelised, unlike the slow progress along the narrow optimisation valley, which is a serial calculation.

\subsection{Combination with Broyden's Method}
As evaluating the Jacobian is more expensive than evaluating the function, several methods exist for estimating the Jacobian based on discovered function values.  Broyden's method \cite{Broyden} updates the approximate Jacobian each iteration via
\[ J_\mathrm{new} = J + \frac{\Delta\f - J\Delta\x}{|\Delta\x|^2}\Delta\x^T, \]
where $\Delta\x$ and $\Delta\f$ are the change in variables and function value, respectively, over the previous step.

\begin{figure}[!htb]
   \centering
   \includegraphics*[width=\columnwidth]{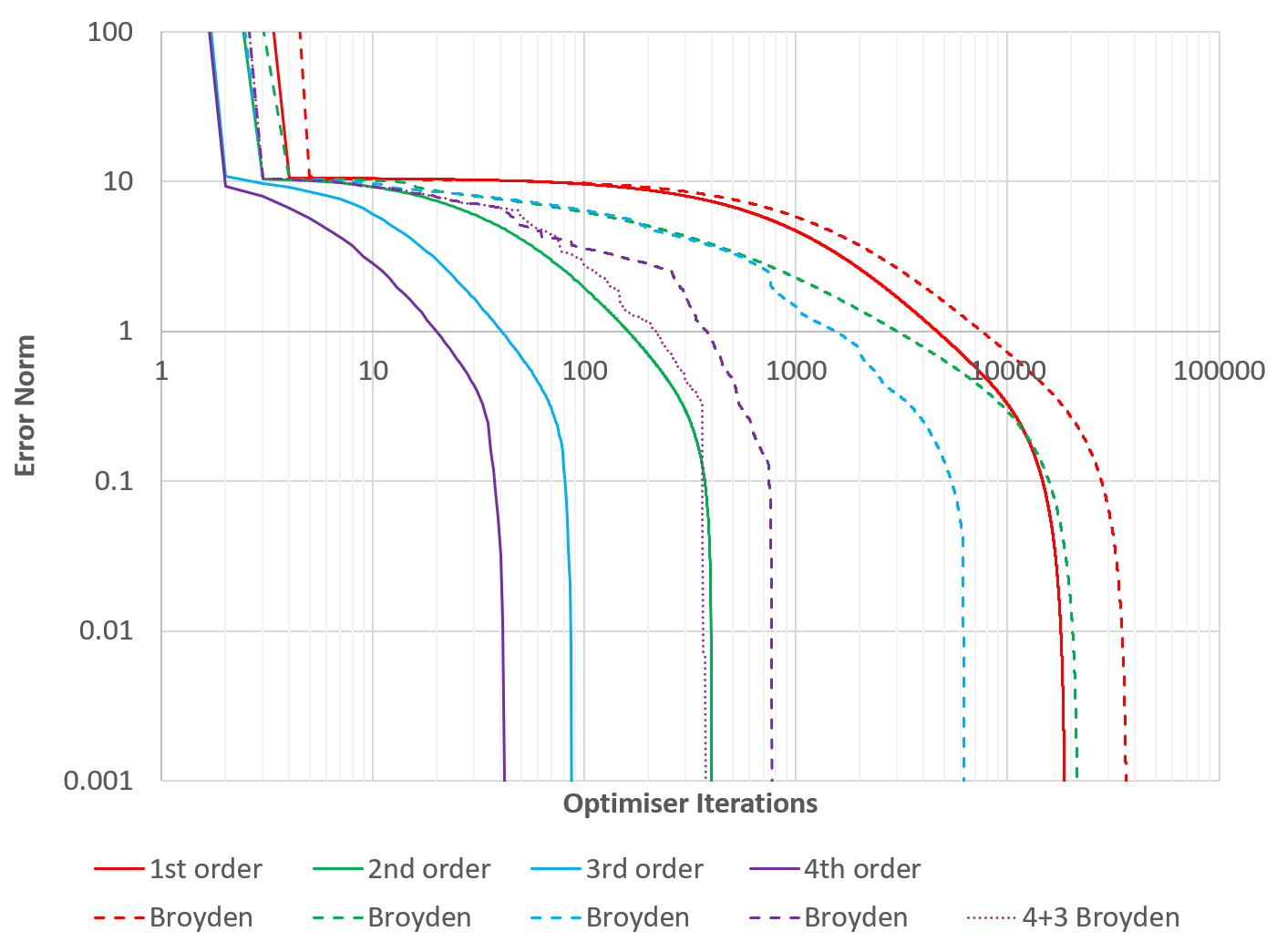}
   \caption{Performance of higher-order corrections on the test problem, with and without using Broyden's method.  `4+3' indicates both the 3$^\mathrm{rd}$ and 4$^\mathrm{th}$ order points were sampled each iteration and the best value taken.}
   \label{fig:testproblem_broyden}
\end{figure}

Figure~\ref{fig:testproblem_broyden} shows the results of Broyden's method on the simple test problem with $K=10^6$.  The solid lines are the same as in Figure~\ref{fig:testproblem}, while the dashed lines only evaluate the Jacobian at the start and use the Broyden update each step thereafter.  The higher-order corrections and Levenburg--Marquardt logic discussed in previous sections are still present, with the approximate Jacobian being used as input rather than the exact one.

\begin{table}[!hbt]
   \centering
   \caption{Test problem convergence times for different order methods using Broyden's method and $K=10^6$.}
   \begin{tabular}{lccc}
       \hline
\textbf{Order} & \textbf{Iterations} & \textbf{Evaluations per iteration} & \textbf{Evaluations} \\
       \hline
1$^\mathrm{st}$ & 36652 & 1 & 36652 \\
2$^\mathrm{nd}$ & 21571 & 2 & 43142 \\
3$^\mathrm{rd}$ & 6211 & 5 & 31055 \\
4$^\mathrm{th}$ & 775 & 9 & 6975 \\
4$^\mathrm{th}$+3$^\mathrm{rd}$ & 376 & 10 & 3760 \\
       \hline
   \end{tabular}
   \label{tab:broyden}
\end{table}

As expected, using an inexact Jacobian results in convergence taking more iterations, but the iterations will be faster by a problem-dependent factor since the Jacobian is not evaluated.  Table~\ref{tab:broyden} compares the convergence times for the Broyden method with different higher-order corrections.  The 4$^\mathrm{th}$ order method significantly reduces not only the number of iterations taken, but the number of function evaluations performed (each 4$^\mathrm{th}$ order step takes 9 evaluations).  A variant of this method that also samples the 3$^\mathrm{rd}$ order point $\f(\x+\cv_1+\cv_2+\cv_3)$, which is not present in the 4$^\mathrm{th}$ order stencil, performs even better.  The best of the two points is taken each iteration.

\section{Performance on a Practical Problem}
These algorithms have also been used on a more complex optimisation problem (which motivated their development).  This problem has 180 parameters and 1200 output variables and features levels of successively more difficult narrow curved valleys in its optimisation space.

The full details of this problem are not the point of this paper but a brief summary will be given here (for physics background, see \cite{HB2023}).  An initial distribution of $N_\mathrm{ions}=400$ Ca$^+$ ions is accelerated through a potential of 1\,kV and given a $\pm$2\% energy chirp.  It is transported through a curved electrostatic channel, where the electric field is produced by 15 rings of 12 configurable electrodes.  Each electrode is modelled as a point charge and these 180 charges are the optimisation variables.  The output vector whose norm should be minimised contains the $(x,y,z)$ position coordinates of the ions on exiting the channel (so $3N_\mathrm{ions}$ entries in all), with the bunch centroid subtracted.  The ions do not interact in this model and their trajectories are calculated by the 4$^\mathrm{th}$ order Runge--Kutta method \cite{OriginalRK,NumRecipesRK} with a timestep $\delta t=10^{-7}$ seconds.

This problem is interesting because focussing the ions to a point in the linear dynamics approximation can be done with standard optics, but optical aberrations at higher order will remain, for example spherical aberration from large angle effects.  These higher order aberrations can also be corrected by careful choice of the electrode voltages, although this gets more difficult the smaller the focal point becomes and the more aberrations have to be cancelled simultaneously.

\begin{figure}[!htb]
   \centering
   \includegraphics*[width=\columnwidth]{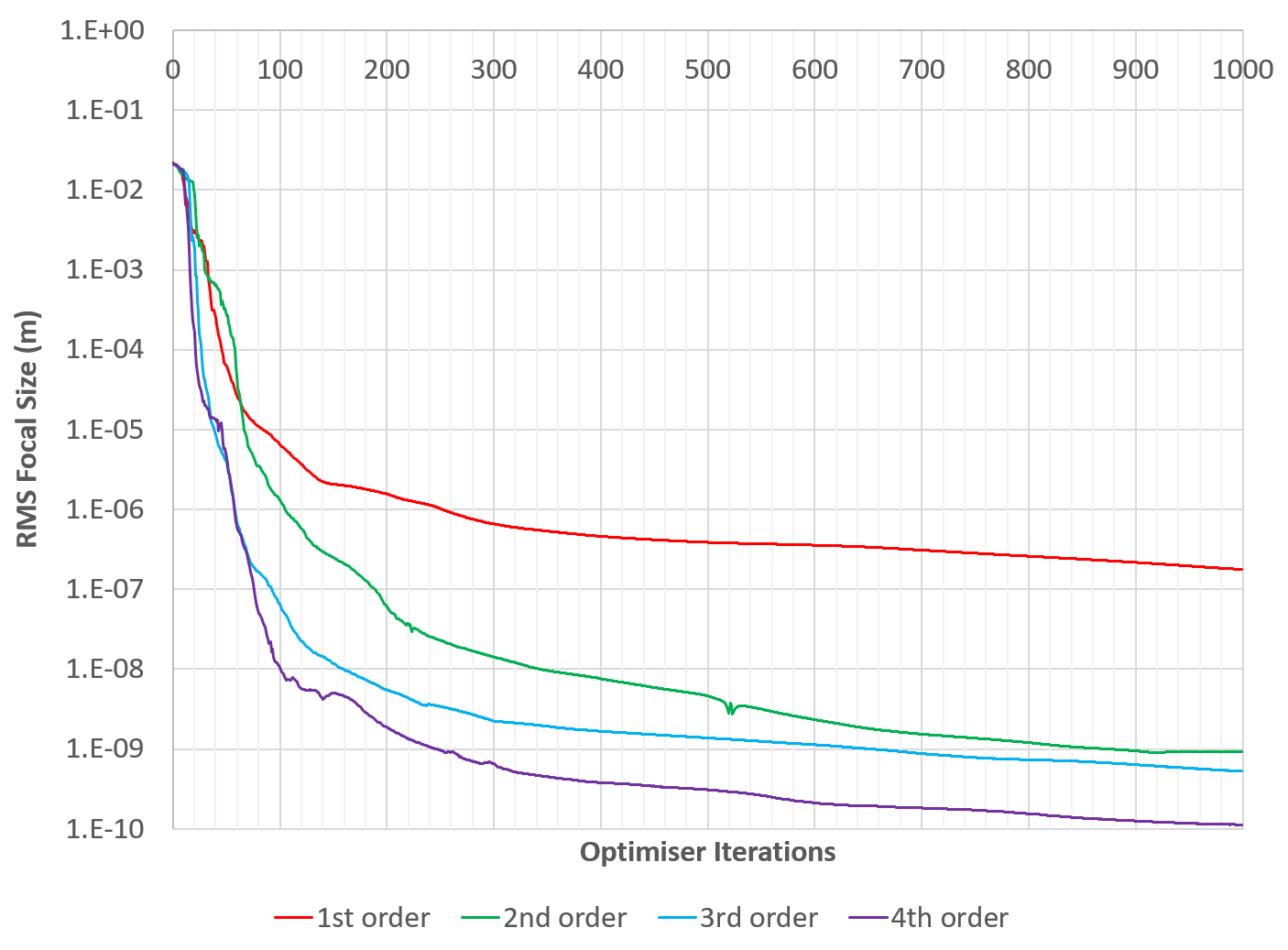}
   \caption{Performance of higher-order corrected Levenburg--Marquardt methods on a physics problem with 180 parameters and 300 output variables.}
   \label{fig:realproblem_iters}
\end{figure}

The figure of merit is the RMS focal size of the ion bunch, which is $\frac1{\sqrt{N_\mathrm{ions}}}$ of the norm $|\f|$.  Figure \ref{fig:realproblem_iters} shows how this is reduced by the Levenburg--Marquardt optimisation method with various levels of higher order correction.  To maintain accuracy, the Jacobian was calculated with the chain rule (through all Runge--Kutta steps) here, rather than finite differences.

Figure \ref{fig:realproblem_iters} shows that higher order methods significantly accelerate the optimisation progress, even on this more complex problem.  All methods eventually slow down here, which is suspected to be because the valley anisotropy continues to increase at lower focal sizes, leading to more steps taken (as in Figure \ref{fig:anisotropy}).

\begin{figure}[!htb]
   \centering
   \includegraphics*[width=\columnwidth]{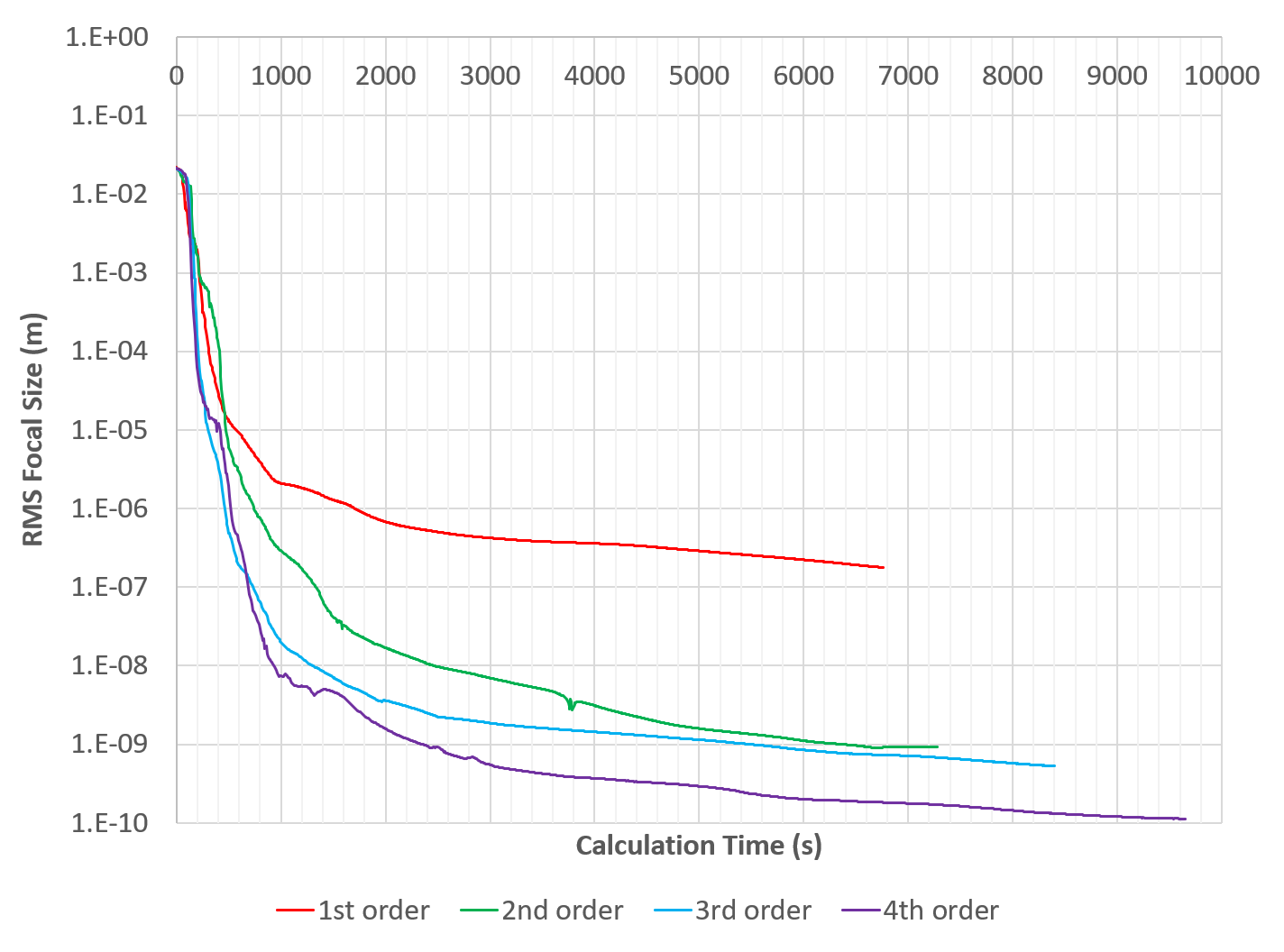}
   \caption{Calculation time vs. performance for higher-order methods on a physics problem with 180 parameters.}
   \label{fig:realproblem_calctime}
\end{figure}

One potential concern with these higher-order methods is that the additional evaluations of $\f$ in the stencil will cost too much time to make the method worth using.  To measure this effect, the optimised focal size is plotted as a function of calculation time in Figure \ref{fig:realproblem_calctime}.  The differences in total execution time can be seen at the right-hand end of each line, which is after 1000 iterations.  The fourth order method takes 44\% longer per step but still manages to pull ahead of the other methods in real time for focal sizes below $10^{-7}$\,m.  Table \ref{tab:realproblem_calctime} gives the amount of time for each method to reach certain RMS focal sizes.

\begin{table}[!hbt]
   \centering
   \caption{Time taken to reach various focal sizes in the ion focussing problem with 180 parameters.}
   \begin{tabular}{lcccc}
       \hline
\textbf{Order} & \textbf{Iterations} \\
 & \textbf{to $<10^{-6}$\,m} & \textbf{to $<10^{-7}$\,m} & \textbf{to $<10^{-8}$\,m} & \textbf{to $<10^{-9}$\,m} \\
       \hline
1$^\mathrm{st}$ & 252 & $>$1000 & $>$1000 & $>$1000 \\
2$^\mathrm{nd}$ & 105 & 189 & 344 & 877 \\
3$^\mathrm{rd}$ & 58 & 93 & 160 & 654 \\
4$^\mathrm{th}$ & 58 & 76 & 101 & 247 \\
       \hline
\textbf{Order} & \textbf{Calculation Time (s)} \\
 & \textbf{to $<10^{-6}$\,m} & \textbf{to $<10^{-7}$\,m} & \textbf{to $<10^{-8}$\,m} & \textbf{to $<10^{-9}$\,m} \\
       \hline
1$^\mathrm{st}$ & 1707.285 & --- & --- & --- \\
2$^\mathrm{nd}$ & 742.992 & 1337.178 & 2470.101 & 6384.886 \\
3$^\mathrm{rd}$ & 461.045 & 739.670 & 1308.115 & 5497.770 \\
4$^\mathrm{th}$ & 525.800 & 690.696 & 934.539 & 2358.116 \\
       \hline
   \end{tabular}
   \label{tab:realproblem_calctime}
\end{table}

\subsection{Practical Problem in Combination with Broyden's Method}
The more complex optimisation problem was attempted using Broyden's method, only evaluating the Jacobian on the first step.  This approach was unsuccessful, resulting in a 5\% decrease in figure of merit in the first 5--10 iterations, followed by either coming to a complete halt, or a sequence of tiny steps decreasing by $\sim$$10^{-10}$ each time.  These steps are not long enough for the higher order corrections to have any effect, nor are they likely to improve the Jacobian accuracy much.  It is suspected that the Jacobian estimate became too far away from the real value to be useful.

\begin{figure}[!htb]
   \centering
   \includegraphics*[width=\columnwidth]{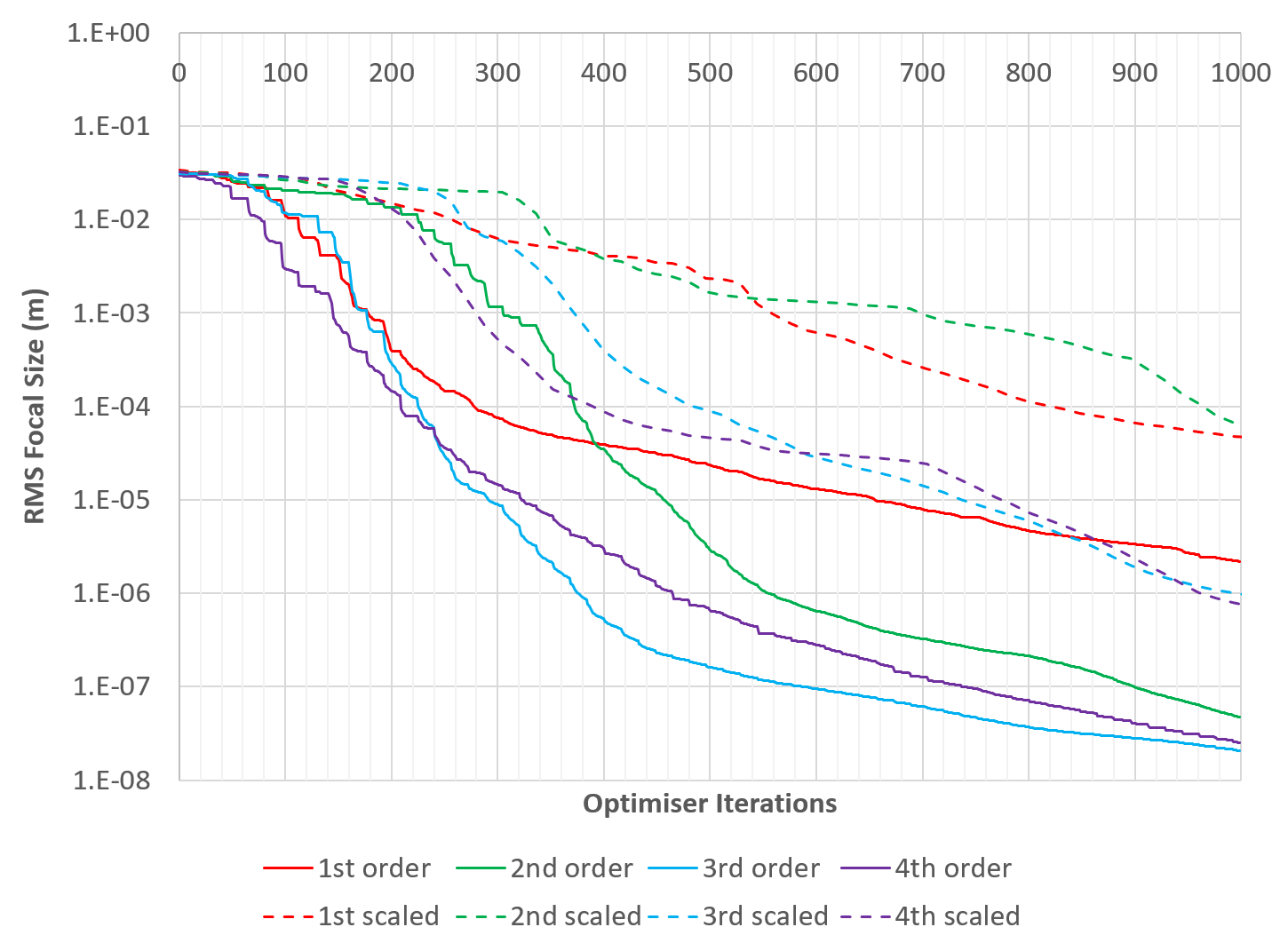}
   \caption{Performance of higher-order corrections on the physics problem, only evaluating the Jacobian every 16$^\mathrm{th}$ iteration and using Broyden's method otherwise.  Dashed lines show the previous results evaluating the Jacobian every time, but with the iteration count multiplied by 16 to make the number of Jacobian evaluations comparable.}
   \label{fig:realproblem_broyden16}
\end{figure}

The optimisation was attempted again with the full Jacobian being evaluated periodically (every 16$^\mathrm{th}$ iteration here) and using Broyden's method for the other steps.  Figure~\ref{fig:realproblem_broyden16} shows that in this case, the intermediate Broyden's method steps do provide some benefit over a direct iteration with the same number of Jacobian evaluations.  It can also be seen that higher-order corrections improve this combined method significantly, although here 3$^\mathrm{rd}$ order wins out slightly over 4$^\mathrm{th}$.  Whether this approach is advantageous will depend on the time taken to evaluate the Jacobian relative to individual function values.  The downward jumps in the graph caused by the Jacobian evaluations suggest they are more important in the earlier stages of the optimisation, where the jumps appear larger relative to the smooth decrease.

\section{Conclusion}
Methods such as Levenburg--Marquardt are not only for curve fitting, but also powerful optimisers where the function to be minimised is a sum of squares.  Like many optimisers, they can get stuck for long periods of time in `curved narrow valleys'.  This paper derives higher-order corrections beyond the already known second order \cite{GeodesicAccelLett,GeodesicAccelE} that further accelerate the optimiser performance in these situations.  A general formula for deriving the $n^\mathrm{th}$ order correction is given, with suggestions on how to build finite difference stencils to evaluate it.

These successful methods have been derived using the concept of a `natural pathway' for the optimisation, which is an ordinary differential equation (ODE) that is meant to follow the valleys.  The form of this ODE is chosen somewhat arbitrarily here but it appears to work well, perhaps because it is a continuous version of the optimiser's path in the limit where step size $\epsilon\rightarrow 0$.  Using a higher-order step method on this ODE thus makes the optimiser behave `as if' it had done a large number of very small steps.  

This link with ODEs also suggests potential future work in applying well-known higher-order methods for ODEs such as Runge--Kutta \cite{OriginalRK,NumRecipesRK} or Bulirsch--Stoer \cite{BulirschStoer} to difficult optimisation problems, as well as methods for stiff ODEs.  In this paper the RK4 method was not preferred because it would have required four evaluations of the Jacobian, which is still much more expensive than the eight additional function points in the stencil of the fourth order method.

\end{document}